\RequirePackage{etex}

\documentclass[10pt]{article}
\usepackage{amsmath, amsthm, amssymb, mathrsfs}
\usepackage{dcpic, pictex, pinlabel, setspace, rotating}
\usepackage{array, verbatim, bbm, faktor}
\usepackage{graphicx}
\usepackage{tabularx}
\usepackage{color}
\usepackage{amscd}
\usepackage{tikz}
\usetikzlibrary{matrix,arrows}
\usepackage[all]{xypic}
\usepackage{xy}

\theoremstyle{definition}

\theoremstyle{remark}

\numberwithin{equation}{section}
\numberwithin{figure}{section}

\def\hline{\bigskip\hrule\bigskip}  

\begin{document}

\title{Group Width}

\author{\sc\small Michael H. Freedman}


\date{}
\maketitle

\begin{abstract}
    \centering
    \begin{minipage}{0.5\textwidth}
        There are many ``minimax" complexity functions in mathematics: width of a tree or a link, Heegaard genus of a $3$-manifold, the Cheeger constant of a Riemannian manifold. We define such a function $w$, ``width", on countable (or finite) groups and show $w(\mathbb{Z}^k)=k-1$.
    \end{minipage}
\end{abstract}

\vspace{1cm}


Let $K$ be a countable (or finite) simplicial complex and give the real line $\mathbb{R}$ the cell structure with the integers $\mathbb{Z}\subset\mathbb{R}$ the vertices. Abusing the usual terminology, we call any simplicial map $f:K\rightarrow\mathbb{R}$ ``morse."\\ \\
\textbf{Definition 1}. \textit{Connected width rank}, $cwr(K):=$ $\rm{min}_{f\,,\rm{morse}}$ $\rm{max}_{i\in\mathbb{Z}}$ rank$(inc_{\#}(\pi_1C))$, where $C$ is some connected component of $f^{-1}[i,i+1]$, and rank means the smallest number of generators of a given group. The inclusion is $C\subset{K}$, and $inc_{\#}(\pi_1C)$ is a subgroup of $\pi_1{K}$. \\ \\
\textbf{Definition 2}. Given a countable group $G$, its \textit{width}, $w(G)$, is the minimum of $cwr(K)$ over all $K$ with $\pi_1(K)\cong{G}$. \\ \\

Clearly free groups (and only free groups) have width zero. Let's work out $w(\mathbb{Z}^k)$, the width of the free abelian group. Consider $f:K\rightarrow\mathbb{R}$ with $\pi_1(K)\cong\mathbb{Z}^k$. Define the quotient (bipartite) graph $Q_f$ by taking an edge for each connected component of $f^{-1}(i)$ and a vertex for each component of $f^{-1}[i,i+1]$, $i\in\mathbb{Z}$, and gluing the latter to the former according to inclusion.

There is an induced map $\theta:K\rightarrow{Q_f}$ and an epimorphism of groups: $$\pi_1K\twoheadrightarrow{\pi_1Q_f},$$ so in our case: $\pi_1K\cong{\mathbb{Z}_{1}^{k}}$, and $Q_f$ is either contractible (a tree) or $Q_f\simeq{S^1}$, is homotopy equivalent to the circle.

First consider the case $Q_f\simeq{p{t}}$. Take a minimal connected subgraph $p\subset{Q_f}$ so that $H_1(\theta^{-1}(p);{Q})\rightarrow{H_1(K,Q)}$ is onto, where $Q$ denotes the rationals. To define ``minimal" we order subgraphs by inclusion. We show that any such $p$ is a single vertex. Suppose that $p$ is minimal but larger than a single vertex. Cut $p$ at the midpoint of some edge $e$ to obtain the complementary subtrees $p_1,p_2\subset{p}$. Let the inverse images under $\theta$ be $P_1,P_2\subset{P}$. Applying the $Q$-homology Mayer-Vietoris sequence to the inclusions (and using connectivity of $P_1\cap{P_2}$), we find that there are classes $b_1\in{H_1(P_1;{Q})}$ and $b_2\in{H_1(P_2;{Q})}$ so that image$(b_1)$ and image$(H_1(P_2;Q))$ are rationally independent in $H_1(K;Q)$ and image$(b_2)$ and image$(H_1(P_1;Q))$ are also independent in $H_1(K;Q)$. Let $\beta_1(\beta_2)\subset{P_1(P_2)}$ be corresponding loops carrying $b_1(b_2)$. The commutation of $\beta_1$ and $\beta_2$ in $\pi_1(K)$ conflicts usefully with the following lemma. Let $T^{+}$ be the $2$-torus $S^1\times{S^1}$ with ``flanges" glued to the factor circles: $$T^+=S^1\times{S^1}\bigcup\limits_{x\times\ast\equiv{x}\times\ast\times{0}}S^{1}\times\ast\times[0,1]\bigcup\limits_{\ast\times{x}\equiv\ast\times{x}\times{0}}\ast\times{S^1}\times[0,1],$$ Denote $S^{1}\times\ast\times{1}$ by $\alpha_1$ and $\ast\times{S^1}\times{1}$ by $\alpha_2$.\\ \\

\pagebreak

\noindent\textbf{Lemma 3}. It is not possible to cover $T^{+}$ by open sets $\mathcal{U}_1$ and $\mathcal{U}_2$ with $\alpha_1\in\mathcal{U}_1$ and $\alpha_2\in\mathcal{U}_2$ with image$(H_1(\mathcal{U}_1,Q))\subset{H_1(T^{+},Q)}$ and image$(H_1(\mathcal{U}_2,Q))\subset{H_1(T^{+},Q)}$, each rank one.\\
\textit{Proof}. This is an exercise in Lusternick-Shirleman category. Consider the cup-product diagram: \vspace{1cm}
\begin*\textcolor{white}{{figure}[htpb]}
\labellist \small\hair 2pt

  \pinlabel $\text{$H^1(T^{+},U_1;Q)$}$ at 118 86
  \pinlabel $\text{$H^1(T^{+};Q)$}$ at 118 31
  \pinlabel $\text{$H^1(T^{+},U_2;Q)$}$ at 300 86
  \pinlabel $\text{$H^1(T^{+};Q)$}$ at 300 31
  \pinlabel $\text{$H^2(T^{+},T^{+};Q)\cong{0}$}$ at 500 86
  \pinlabel $\text{$H^2(T^{+};Q)$}$ at 495 31
  \pinlabel $\text{$\times$}$ at 207 86
  \pinlabel $\text{$\times$}$ at 207 31
  \pinlabel $\text{$\times$}$ at 207 4
  \pinlabel $\text{$\reflectbox{\rotatebox[origin=c]{90}{$\in$}}$}$ at 118 13
  \pinlabel $\text{$\reflectbox{\rotatebox[origin=c]{90}{$\in$}}$}$ at 302 13
  \pinlabel $\text{$\reflectbox{\rotatebox[origin=c]{90}{$\in$}}$}$ at 495 13
  \pinlabel $\text{$\hat{\alpha_1}$}$ at 118 0
  \pinlabel $\text{$\hat{\alpha_2}$}$ at 302 0
  \pinlabel $\text{$1$}$ at 495 0
  \pinlabel $\text{$i_1$}$ at 128 56
  \pinlabel $\text{$i_2$}$ at 315 56
  \pinlabel $\text{$\reflectbox{\rotatebox[origin=c]{270}{$\xrightarrow{\hspace*{5mm}}$}}$}$ at 115 56
  \pinlabel $\text{$\reflectbox{\rotatebox[origin=c]{270}{$\xrightarrow{\hspace*{5mm}}$}}$}$ at 302 56
  \pinlabel $\text{$\reflectbox{\rotatebox[origin=c]{270}{$\xrightarrow{\hspace*{5mm}}$}}$}$ at 495 56
  \pinlabel $\text{$\xrightarrow{\hspace*{2cm}}$}$ at 395 86
  \pinlabel $\text{$\xrightarrow{\hspace*{2cm}}$}$ at 395 31
  \pinlabel $\text{$\small{|}\normalsize\hspace{-1.8mm}\xrightarrow{\hspace*{2cm}}$}$ at 395 0

\endlabellist
\centering
\includegraphics[scale=.8]{background}
\end*\textcolor{white}{{figure}}


\noindent Using the exact sequence of pairs, Image$(i_1)=$ Span$(\hat{\alpha_1})$ and Image$(i_2)=$ Span$(\hat{\alpha_2})$, where $\hat{
\alpha_i}$ denotes the Poincar\'{e} dual to the loops $\alpha_i$, $i=1,2$. The factoring of the cup product $\hat{\alpha_1}\smile\hat{\alpha_2}=1$ through zero is a contradiction. \qed \\ \\

Since $\pi_1(K)$ is abelian, there is a map $g:T^{+}\rightarrow{K}$ carrying $\alpha_1$ to $\beta_1$ and $\alpha_2$ to $\beta_2$. Taking $\theta^{-1}$, cutting the edge $e$ divides $K$ into $K_1$ and $K_2$, containing $P_1$ and $P_2$ (resp.). Let $K_1^{+}$ and $K_2^{+}$ be homotopy equivalent open sets containing $K_1$ and $K_2$ (resp.), $K_1^{+}\simeq{K_1}$ and $K_2^{+}\simeq{K_2}$. Setting $\mathcal{U}_i=g^{-1}(K_i^{+})$, $i=1,2$, contradicts the lemma, showing $P$ consists of a single vertex $v$. Thus if $Q_f$ is a tree, a finite index sublattice $L$ of $\pi_1(K)\cong\mathbb{Z}^{k}$ must be generated by $\theta^{-1}(v)$, and rank$(L)=k$.

Next consider the case $Q_f\simeq{S^{1}}$. Again take a minimal $p\subset{Q_f}$. $p$ must contain the essential loop $\gamma\subset{Q_f}$, otherwise the epimorphism $H_1(K,Q)\overset{\theta_{\ast}}\twoheadrightarrow{H_1(Q_f;Q)}$ would factor through a trivial $H_1(p;Q)$. By the preceeding argument, $p=\gamma$; we may trim off leaves of $Q_f$ by arguing they cannot increase the image in the rationalized fundamental group, $H_1(K;{Q})$.

Let $v_n=v_0,v_1,...,v_{n-1}$ be the vertices on $p$ and $V_n=V_0,V_1,...,V_{n-1}$ the $\theta$-preimages. We claim that for $0\leq{i,j}\leq{n-1}$, Image$(H_1(V_i;\mathbb{Z}))=$ image$(H_1(V_j;\mathbb{Z}))\subset{H_1(K)}$. To see this, note that for any loop $\delta\subset{V_i}$, there is a map of a torus $h:S^1\times{S^1}\rightarrow{K}$ with $\theta{h}(S^{1}\times\ast)$ parameterizing $\gamma$ and $h(\ast\times{S^1})$ parameterizing $\delta$. Using transversality, we may arrange that corresponding to the center point $\widehat{e_1},...,\widehat{e_n}$ of each edge in $p$, $h^{-1}(\widehat{e_k})$ is a $1$-manifold in $S^1\times{S^1}$ meeting $S^1\times\ast$ transversely in a single point. These $1$-manifolds all (up to sign) represent the same class $inc_{\ast}[\delta]\in{H_1(K;\mathbb{Z})}$ since they are homologous on $S^1\times{S^1}$.

Using the connectivity of $\theta^{-1}(v_k)$, $k=0,..,{n-1}$, and again, the Mayer-Vietoris sequence, we see that $inc_{\ast}{H_1(V_0\cup\tilde{\gamma};{Q})}=H_1(K;{Q})$ where $\tilde{\gamma}$ is some lift of $\gamma$, $\theta(\tilde{\gamma})=\gamma$. Similarly, for all $V_k$, $1\leq{k}\leq{n-1}$. It is also clear that $inc_{\ast}[\tilde{\gamma}]$ and $inc_{\ast}H_1(V_0;{Q})$ must be indepedent in $H_1(K;{Q})$, otherwise a homotopy, in $K$, of a multiple of $\tilde{\gamma}$ into $V_0$, would, under $\theta$, constitute a null homotopy in $Q_f$ of a multiple of the essential cycle $\gamma$. Thus $inc_{\ast}H_1(V_0;{Q})$ has rank$=k-1$ in $H_1(K;{Q})$.

We have shown that if $Q_f\simeq\ast$, then some component $C$ carries all of $\pi_1(K)\otimes{Q}\cong{H_1(K;{Q})}$ and if $Q_f\simeq{S^1}$ then some component $C$ carries a rank $k-1$ subgroup. Since the obvious Morse function on the $k$-torus $T^{k}$ has all levels (even critical levels) carrying a rank $k-1$ subspace of $H_1(T^{k},\mathbb{Z})$ (and has $Q_f\equiv{S^1}$), we conclude that $w(\mathbb{Z}^k)=k-1$. 

\section*{Extensions}

\indent \indent For a finite abelian group $A$, $w(A)=$ rank$(A)$. The proof is similar to the computation of $w(\mathbb{Z}^k)$ except for two modifications. First, $Q_f$ is now certainly a tree so only that case requires generalization. Second, in all computations, the rationals $Q$ should be replaced with the field $\mathbb{Z}/q\mathbb{Z}$, where $q$ is a prime contained in the factorization of order$(A)$ at least as often as any other prime.

Combining the arguments for both free and torsion cases, one finds that for a finitely generated, but infinite, abelian group $B$, that $w(B)=$ rank$(B)-1$.

It is easy to say a little more about finite groups: a finite group $F$ of width one is cyclic. To prove this, assume $\pi_1(K)\cong{F}$ and $f:K\rightarrow\mathbb{R}$ exhibits the width of $F$ to be one. Choose a maximal subtree $p\subset{Q_f}$ with respect to the property that $\pi_1(P)$ has cyclic image $X$ in $\pi_1(K)$, where $P=\theta^{-1}(p)$. Let $p^{+}$ be $p$ union an adjacent $1$-simplex $e$ of $Q_f$ and let $P^{+}=\theta^{-1}(p^{+})$. Write $P^{+}=P\cup{C}$, where $C=\theta^{-1}(e)$. Note image$(\pi_1(P^+))=H\subset\pi_1(K)$ is not cyclic, but image$(\pi_1(C))=:Y\subset{F}$ is cyclic. Let ${Z}\subset{F}$ be the cyclic group ${Z}=X\cap{Y}\subset{F}$ and let $G:=X\underset{{Z}}\ast{Y}$ be the abstract free product with amalgamation. There is an epimorphism $\gamma:G\twoheadrightarrow{H}$. Since $G$ is infinite and $H$ is finite, there must be a nontrivial relation $R\in\ker{\gamma}$. Since ${Z}\cap\ker{\gamma}=\{id.\}$, $R$ can be written as a cyclically reduced word alternating ``letters" from $X\setminus{Z}$ and $Y\setminus{Z}$. Think of $R$ as a map $R:D^2\rightarrow{K}$, which on the boundary maps to a wedge of circles $S^1\bigvee{S^1}$, the first summand lying in $P$ and the second summand in $C$. Make $R$ transverse to $P\cap{C}$ and consider an innermost arc $\omega\subset{D^2}$, $\omega\subset{R^{-1}(P\cap{C})}$. The subdisk $\Delta\subset{D^2}$ between $\omega$ and $\partial{D^2}$ determines a (pointed) homotopy of some letter of $R$ into $P\cap{C}$. Since image$(\pi_1(P\cap{C}))\subset{Z}\subset{F}$, this contradicts the form of $R$, i.e. that its letters lie in $(X\setminus{Z})\amalg(Y\setminus{Z})$. It follows that $p=Q_f$ and $F=\pi_1(K)$ is cyclic.

Formal properties of width include: $$w(G_1\times{G_2})\leq{w(G_1)+}\text{ rank}(G_2)$$ and $$w(G_1\ast{G_2})=\text{ max}\{w(G_1),w(G_2)\}.$$

To prove the latter, given $f:K\rightarrow\mathbb{R}$ with $\pi_1(K)\cong{G_1\ast{G_2}}$, one may precompose with the covering $\delta_j:K_j\rightarrow{K}$, $\pi_1(K_j)\cong{G_j}$, to obtain $f_j=f\circ\delta_j$, $j=1$ or $2$. By Grushko's decomposition theorem, for any connected component $C_j$ of $f_{j}^{-1}[i,i+1]$, the image $H_j$ of $\pi_1(C_j)$ in $\pi_1(K)$ is a free summand of the corresponding image $H$ of $\pi_1(C)$ in $\pi_1(K)$, where $C=\delta_j(C_j)$. Consequently, rank$(H_j)\leq$ rank$(H)$, establishing $w(G_1\ast{G_2})\geq{\text{max}}\{w(G_1),w(G_2)\}$. The opposite inequality is immediate.

\section*{Applications}

\indent\indent The computation $w(\mathbb{Z}^k)=k-1$ immediately gives negative answers to two MathOverflow questions: mathoverflow.net/questions/30567/ and mathoverflow.net/questions/42629/. More specifically, for dimension $d\geq{4}$ consider a smooth closed $d$-manifold $M$ with $\pi_1(M)=\mathbb{Z}^k$. Any Morse function on $M$ must have some connected component $C$ of some level with first Betti number $b_1(C)\geq{k-1}$. It requires only a little thought to see that this estimate also applies to generic levels. The ``complexity" of connected levels is thus seen to increase with $k$. 

Secondly, if $M$ is divided up into connected ``blocks" along codimension$=1$ manifold faces, at least one block must have $b_1($block$)\geq{k-1}$ (blocks will map to vertices of $Q_f$, their faces to edges of $Q_f$). Product collars can be added along faces to build a simplicial Morse function $M\rightarrow{\mathbb{R}}$ as in Definition 1, with all blocks corresponding to components $C$. Thus general $d$ manifolds cannot be cut into simple pieces, comprising only a finite number $n_d$ of diffeomorphism types, with purely $(d-1)$-manifold cuts, as was asked. If the (second) question, instead, permitted gluing along codimension one \textit{and} two faces, it would not be touched by this group theoretic method and appears open. Also, restricting the question to simply connected manifolds would require a different method. This question looks difficult in the case of simply connected (smooth) $4$-manifolds.

\section*{Acknowledgements}
I would like to thank Ian Agol for discussions and the references to MathOverflow. \\ \\ \\





\noindent
\noindent\small\sc{Microsoft Station Q, University of California, Santa Barbara, CA 93106} \\ \rm
\textit{Email address:}\texttt{ michaelf@microsoft.com}

\end{document}